\documentclass[11pt]{article}

\usepackage[a4paper,margin=1in]{geometry}
\usepackage{amsmath,amssymb,amsthm}
\usepackage{booktabs}
\usepackage{longtable}
\usepackage{array}
\usepackage{hyperref}
\usepackage{enumitem}

\title{A Distinct Covering System with Minimum Modulus $7$ and Minimal Least Common Multiple $10080$}

\author{
\begin{tabular}{c}
Shiliang Zhang\textsuperscript{1}
\quad
Jiheng Zhang\textsuperscript{1}\\[0.4em]
\textsuperscript{1}Institute of Mathematical Sciences, ShanghaiTech University
\end{tabular}
}
\date{\today}

\newcommand{\ContactInfo}{{
\bigskip\footnotesize

\medskip
\noindent S.~Zhang,
\textsc{Institute of Mathematical Sciences\\
ShanghaiTech University\\
No.393 Middle Huaxia Road, Pudong New District, 
Shanghai, China}\par\nopagebreak
\noindent\textsc{Email:} \texttt{zhangshl2023@shanghaitech.edu.cn}

\bigskip
\noindent J.~Zhang,
\textsc{Institute of Mathematical Sciences\\
ShanghaiTech University\\
No.393 Middle Huaxia Road, Pudong New District, 
Shanghai, China}\par\nopagebreak
\noindent\textsc{Email:} \texttt{zhangjh2023@shanghaitech.edu.cn}

}}

\newtheorem{theorem}{Theorem}
\newtheorem{lemma}{Lemma}

\newtheorem{definition}{Definition}
\newtheorem{example}{Example}

\newcommand{\Z}{\mathbb{Z}}
\newcommand{\lcm}{\operatorname{lcm}}

\begin{document}

\maketitle

\begin{abstract}
We determine the minimum possible least common multiple of a distinct covering system whose minimum modulus is \(7\). 
Klein previously constructed such a system with least common multiple \(15120\) and conjectured that this value was minimal.
We give a construction with least common multiple \(10080\), and we prove that no smaller least common multiple can occur. 
The proof is organized as a successive filtering argument. 
Starting from the possible multiples of \(7\) below \(10080\), we first apply a reciprocal-sum filter, then a divisor-completed integer-programming filter, then a stronger partial-sum filter. 
The few remaining hard cases are finally certified by complete Gurobi computations.
\end{abstract}

\section{Introduction}\label{sec:introduction}

Covering systems were introduced by Erd\H{o}s in 1950 in connection with Romanov's question on integers of the form \(2^k+p\), where \(p\) is prime \cite{Erdos1950}. 
Since then, they have become a central object in the study of arithmetic progressions covering the integers. 
For general background and surveys, see \cite{Balister2024,PorubskySchonheim2002}. 
One of the main themes in the area is to understand how large the least modulus in a distinct covering
system can be. 
Hough resolved Erd\H{o}s's minimum-modulus problem by proving that the least modulus in any distinct covering system is bounded above by an absolute constant; more precisely, it is less than \(10^{16}\) \cite{Hough2015}.

In parallel with this qualitative boundedness question, many authors have constructed explicit distinct covering systems with prescribed least modulus.
Erd\H{o}s gave an early construction with least modulus \(3\) \cite{Erdos1950}.
Swift later gave constructions with least modulus \(4\) and \(6\) \cite{Swift1954}. 
Further examples were obtained by Churchhouse, Krukenberg, Choi, Morikawa, Gibson, Nielsen, and Owens, culminating in the current record construction with least modulus \(42\) \cite{Churchhouse1968,Krukenberg1971,Choi1971,Morikawa1981,Gibson2009,Nielsen2009,Owens2014}.

A more refined problem fixes the minimum modulus and asks how small the least common multiple of all moduli can be. 
For a positive integer \(m\), let \[ L_{\min}(m) \] denote the smallest possible least common multiple among all distinct covering systems whose minimum modulus is \(m\). 
Dalton and Trifonov proved \[ L_{\min}(3)=120, L_{\min}(4)=360, \] and Klein later proved \[ L_{\min}(5)=1440, L_{\min}(6)=5040. \]
Thus the least-common-multiple problem is known for minimum modulus at most \(6\) \cite{DaltonTrifonov2022,Klein2026}.

The next open case is \(m=7\). 
Churchhouse's earlier work included a system with minimum modulus \(7\), but with least common multiple \(30240\) \cite{Churchhouse1968}. 
Klein improved this to a distinct covering system with minimum modulus \(7\) and least common multiple \(15120\), and conjectured that \(15120\) was minimal \cite[Conjecture~2]{Klein2026}.

The purpose of this paper is to show that the conjectured value \(15120\) is not minimal. 
We construct a distinct covering system with minimum modulus \(7\) and least common multiple \(10080\), and we prove that every distinct covering
system with minimum modulus \(7\) has least common multiple at least \(10080\).
Equivalently, we prove the following theorem.

\begin{theorem}\label{thm:main}
The smallest possible value of the least common multiple of a distinct covering system with minimum modulus \(7\) is \[ L_{\min}(7)=10080. \]
\end{theorem}

The proof is organized around a filtering strategy.
Since the modulus \(7\) itself must occur, any possible least common multiple must be divisible by \(7\).
The known value \(L_{\min}(6)=5040\), together with a preliminary reciprocal-sum argument, places the search above \(5040\). The construction in this paper gives the upper bound \(10080\). 
We therefore examine multiples of \(7\) in the range
\[ 5040\leq L<10080. \]
The first filter is a reciprocal-sum filter, which eliminates all but a short explicit list of candidates. 
The second filter is a divisor-completed integer programming formulation, with translation and residue-layer normalizations. 
A third, stronger partial-sum filter rules out two difficult candidates by optimizing the actual covered proportion of a selected subset of moduli. 
Finally, the remaining hard candidates are settled by complete Gurobi computations using the full feasibility test.

\section{Preliminaries}\label{sec:preliminaries}

\begin{definition}
A \emph{covering system} is a finite collection of arithmetic progressions \[ a_i \pmod {m_i},\qquad 1\leq i\leq n, \] where \(m_i\geq 2\) and \(0\leq a_i<m_i\), such that
\[ \bigcup_{i=1}^n \{x\in\Z:x\equiv a_i\pmod {m_i}\} = \Z. \]
It is called \emph{distinct} if the moduli \(m_i\) are pairwise distinct. 
In that case, after reordering, we write \[ 2\leq m_1<m_2<\cdots<m_n. \]
The modulus \(m_1\) is called the \emph{minimum modulus} of the system, and \[ L=\lcm(m_1,\ldots,m_n) \] is the least common multiple of its moduli.
\end{definition}

For a fixed least common multiple \(L\), the covering condition may be checked modulo \(L\): a collection of congruence classes covers \(\Z\) if and only if it covers every residue class in \(\Z/L\Z\). 
Moreover, every modulus appearing in the system must divide \(L\).

We shall repeatedly use the following divisor-completion principle. 
Suppose that a distinct covering system has minimum modulus \(m\) and least common multiple \(L\). 
Then its set of moduli is contained in \[ D_m(L):=\{d:d\mid L,\ d\geq m\}. \]
If a divisor \(d\in D_m(L)\) is missing from the system, then adding one arbitrary congruence class modulo \(d\) preserves the covering property, preserves distinctness, and does not change the least common multiple. 
Therefore, to prove that no system exists for a fixed \(L\), it is enough to prove infeasibility for the divisor-completed verification using all moduli in \(D_m(L)\).

We also use the following reciprocal obstruction, which is often attributed to the Mirsky--Newman theory and is stated in the form needed here in \cite[Lemma~5]{Klein2026}.

\begin{lemma}[Reciprocal obstruction]\label{lem:reciprocal}
If \[ C=\{a_i\pmod {m_i}:1\leq i\leq n\} \] is a distinct covering system, then \[ \sum_{i=1}^n \frac{1}{m_i}>1. \]
\end{lemma}

For the present paper, the known \(m=6\) case is the lower anchor for the search.

\begin{theorem}[Klein {\cite[Theorem~4]{Klein2026}}]\label{thm:klein-six}
The least possible least common multiple of a distinct covering system with minimum modulus \(6\) is \[ L_{\min}(6)=5040. \]
\end{theorem}

We now explain why no value \(L<5040\) needs to be considered for minimum modulus \(7\). 
If such a system existed, then \(7\mid L\). 
By Lemma~\ref{lem:reciprocal} and divisor completion, it would also satisfy \[ \sum_{\substack{d\mid L\\ d\geq 7}}\frac1d>1. \]
A direct enumeration of all \(L<5040\) satisfying these two necessary conditions gives \[ \{1260,1680,2520,3024,3360,3780,4032,4200,4620\}. \]
Every value in this list is divisible by \(6\). 
If a distinct covering system with minimum modulus \(7\) and such an \(L\) existed, then adjoining one additional class, for instance \(0\pmod 6\), would produce a distinct covering system with minimum modulus \(6\) and least common multiple still equal to \(L<5040\), contradicting Theorem~\ref{thm:klein-six}. 
Hence any distinct covering system with minimum modulus \(7\) has \(L\geq 5040\).

\section[Candidate Moduli for m1=7]{Candidate Moduli for \(m_1=7\)}\label{sec:candidates}

The construction in Section~\ref{sec:gurobi} gives a distinct covering system with minimum modulus \(7\) and least common multiple \(10080\). 
Thus, in order to prove Theorem~\ref{thm:main}, it remains to rule out all possible least common multiples \(L<10080\).

By the last paragraph of Section~\ref{sec:preliminaries}, no such \(L\) can be smaller than \(5040\). 
Since the modulus \(7\) occurs in every system whose minimum modulus is \(7\), every possible least common multiple must be divisible by \(7\). 
The raw candidate set is therefore \[ \mathcal R_7 = \{L:5040\leq L<10080,\ 7\mid L\}. \]
The subsequent sections apply increasingly strong filters to \(\mathcal R_7\).
The goal is to prove that every \(L\in\mathcal R_7\) is impossible, while
\(L=10080\) is possible.

\section{The Sum Filter}\label{sec:sum-filter}

The first filter is the \emph{sum filter}. It is the computational use of the reciprocal obstruction from Lemma~\ref{lem:reciprocal}. 
For a proposed minimum modulus \(m\) and a proposed least common multiple \(L\), define \[ S_m(L) = \sum_{\substack{d\mid L\\ d\geq m}}\frac1d. \]
The quantity \(S_m(L)\) is the reciprocal sum over all moduli that could appear after divisor completion. 
If a distinct covering system with minimum modulus \(m\) and least common multiple \(L\) exists, then its moduli form a subset of \(D_m(L)\).
Lemma~\ref{lem:reciprocal} therefore gives the necessary condition \[ S_m(L)>1. \]
Consequently, if \(S_m(L)\leq 1\), then \(L\) is impossible. 
This is called the sum filter since it discards candidates using only the reciprocal sum of the
available divisors.

We apply this filter to the raw candidate set \(\mathcal R_7\). 
Thus \(L\) is retained only if \[ S_7(L) = \sum_{\substack{d\mid L\\ d\geq 7}}\frac1d >1. \]
A direct enumeration gives the following surviving values:
\[ \begin{aligned} \mathcal L_7=\{& 5040,5544,5880,6048,6300,6552,6720,6930,7056,\\ &7392, 7560,8064,8400,8568,8820,9072,9240,9576 \}.
\end{aligned}
\]
Equivalently, the sum filter excludes precisely the values in \(\mathcal R_7\setminus\mathcal L_7\). 
It remains to rule out the eighteen candidates in \(\mathcal L_7\). 
For comparison and for the final construction, we also solve the corresponding verification at \(L=10080\).

\section{The Integer Programming Filter}\label{sec:ip-filter}

The second filter is an integer programming feasibility formulation. 
Its purpose is to test, for a fixed \(L\), whether there is a choice of one residue class for each divisor-completed modulus that covers every residue class modulo \(L\). 
If this formulation is infeasible, then no distinct covering system with minimum modulus \(7\) and least common multiple \(L\) can exist.

For \(L\in\mathcal L_7\cup\{10080\}\), set \[ M_L:=D_7(L)=\{d:d\mid L,\ d\geq 7\}. \]
For each \(m\in M_L\) and each residue \(a\in\{0,1,\ldots,m-1\}\), introduce a binary variable \[ x_{m,a}\in\{0,1\}. \]
The interpretation is that \(x_{m,a}=1\) if and only if the congruence class \(a\pmod m\) is selected.

The integer programming feasibility formulation consists of the following constraints:
\begin{align}
    \sum_{a=0}^{m-1}x_{m,a}&=1,
     m\in M_L, \label{eq:one-residue}\\
    \sum_{m\in M_L}x_{m,\,b\bmod m}&\geq 1,
    b=0,1,\ldots,L-1. \label{eq:cover-residue}
\end{align}
The constraints \eqref{eq:one-residue} enforce that exactly one residue class is chosen for each modulus. 
The constraints \eqref{eq:cover-residue} enforce that each residue \(b\in\Z/L\Z\) is covered by at least one chosen congruence class.
There is no objective function: feasibility gives a divisor-completed covering system, while infeasibility rules out \(L\).

In practice, symmetry reduction is important. 
We use two elementary normalizations. 
The first is translation invariance.

\begin{lemma}[Translation invariance {\cite[Lemma~3]{Klein2026}}]\label{lem:translation}
Let \[  C=\{a_i\pmod {m_i}:1\leq i\leq n\} \]
be a covering system. 
For every integer \(t\), the translated system \[ C+t:=\{a_i+t\pmod {m_i}:1\leq i\leq n\} \] is again a covering system with the same set of moduli.
\end{lemma}

Thus, for pairwise coprime moduli, the Chinese Remainder Theorem allows
simultaneous residue normalization. For instance, whenever \(7,8,9\mid L\), we
may assume that the selected classes include
\[
    6\pmod 7,\qquad 7\pmod 8,\qquad 8\pmod 9.
\]

The second normalization is a residue-layer symmetry.

\begin{lemma}[Residue-layer swap {\cite[Lemma~4]{Klein2026}}]\label{lem:layer-swap}
Let \(C\) be a covering system with least common multiple \(L\), and let \(p\)
be a prime divisor of \(L\). For \(\alpha\in\{0,1,\ldots,p-1\}\), define the
\(p\)-residue layer
\[
    C_p(\alpha)
    :=
    \{a\pmod m\in C:p\mid m,\ a\equiv \alpha\pmod p\}.
\]
Let \(0\leq \alpha_1<\alpha_2\leq p-1\), and let \(t\in[0,L)\) be the unique
integer satisfying
\[
    t\equiv \alpha_2-\alpha_1\pmod {p^{\nu_p(L)}},
    \qquad
    t\equiv 0\pmod {q^{\nu_q(L)}}\quad(q\mid L,\ q\neq p).
\]
Here \(\nu_r(L)\) denotes the exponent of the prime \(r\) in the prime
factorization of \(L\). Then
\[
    \bigl(C\setminus (C_p(\alpha_1)\cup C_p(\alpha_2))\bigr)
    \cup
    (C_p(\alpha_1)+t)
    \cup
    (C_p(\alpha_2)-t)
\]
is again a covering system. In particular, two residue layers modulo \(p\) may
be swapped while preserving the covering property and the set of moduli.
\end{lemma}

The following example illustrates how these normalizations give fixed classes in
the integer programs.

\begin{example}\label{ex:normalization-5040}
Consider \(L=5040=2^4\cdot 3^2\cdot 5\cdot 7\). By translation and CRT
normalization, we may assume that the system contains
\[
    6\pmod 7,\qquad 7\pmod 8,\qquad 8\pmod 9.
\]
Since we work in the divisor-completed setting, the modulus \(35\) occurs. Write
its selected class as \(a\pmod {35}\).

First normalize the residue modulo \(5\). Let \(r\equiv a\pmod 5\). If
\(r\neq 3\), apply Lemma~\ref{lem:layer-swap} with \(p=5\), swapping the
\(5\)-layers \(r\pmod 5\) and \(3\pmod 5\). The required translation satisfies
\[
    t_5\equiv 3-r\pmod 5,\qquad
    t_5\equiv 0\pmod {2^4},\qquad
    t_5\equiv 0\pmod {3^2},\qquad
    t_5\equiv 0\pmod 7.
\]
This operation preserves the already fixed classes modulo \(7\), \(8\), and
\(9\). Hence we may assume \(a\equiv 3\pmod 5\).

Next normalize the residue modulo \(7\). Let \(s\equiv a\pmod 7\). If \(s=6\),
then \(a\pmod {35}\subseteq 6\pmod 7\), so the class modulo \(35\) is redundant;
replacing it by \(33\pmod {35}\) cannot destroy the covering property. If
\(s\neq 6\), then we may use Lemma~\ref{lem:layer-swap} with \(p=7\). If
\(s\neq 5\), swap the \(7\)-layers \(s\pmod 7\) and \(5\pmod 7\), choosing
\[
    t_7\equiv 5-s\pmod 7,\qquad
    t_7\equiv 0\pmod {2^4},\qquad
    t_7\equiv 0\pmod {3^2},\qquad
    t_7\equiv 0\pmod 5.
\]
This does not change the condition \(a\equiv 3\pmod 5\), and it does not affect
the fixed class \(6\pmod 7\). Therefore we may assume
\[
    a\equiv 3\pmod 5,\qquad a\equiv 5\pmod 7.
\]
By the Chinese Remainder Theorem, this is equivalent to
\[
    a\equiv 33\pmod {35}.
\]
Thus the preset \(33\pmod {35}\) is valid for the test of \(L=5040\).
\end{example}

Applying the same normalization procedure to the candidates gives the following
fixed congruence classes.

\begin{longtable}{c l}
\toprule
\(L\) & Fixed congruence classes \\
\midrule
5040  & \(6\pmod 7,\,7\pmod 8,\,8\pmod 9,\,33\pmod {35}\) \\
5544  & \(6\pmod 7,\,7\pmod 8,\,8\pmod 9,\,10\pmod {11}\) \\
5880  & \(6\pmod 7,\,7\pmod 8,\,14\pmod {15}\) \\
6048  & \(6\pmod 7,\,7\pmod 8,\,8\pmod 9\) \\
6300  & \(6\pmod 7,\,8\pmod 9,\,9\pmod {10}\) \\
6552  & \(6\pmod 7,\,7\pmod 8,\,8\pmod 9,\,12\pmod {13}\) \\
6720  & \(6\pmod 7,\,7\pmod 8,\,14\pmod {15}\) \\
6930  & \(6\pmod 7,\,8\pmod 9,\,9\pmod {10},\,10\pmod {11}\) \\
7056  & \(6\pmod 7,\,7\pmod 8,\,8\pmod 9\) \\
7392  & \(6\pmod 7,\,7\pmod 8,\,10\pmod {11}\) \\
7560  & \(6\pmod 7,\,7\pmod 8,\,8\pmod 9,\,33\pmod {35}\) \\
8064  & \(6\pmod 7,\,7\pmod 8,\,8\pmod 9\) \\
8400  & \(6\pmod 7,\,7\pmod 8,\,14\pmod {15}\) \\
8568  & \(6\pmod 7,\,7\pmod 8,\,8\pmod 9,\,16\pmod {17}\) \\
8820  & \(6\pmod 7,\,8\pmod 9,\,9\pmod {10}\) \\
9072  & \(6\pmod 7,\,7\pmod 8,\,8\pmod 9\) \\
9240  & \(6\pmod 7,\,7\pmod 8,\,10\pmod {11},\,14\pmod {15},\,33\pmod {35}\) \\
9576  & \(6\pmod 7,\,7\pmod 8,\,8\pmod 9,\,18\pmod {19}\) \\
10080 & \(6\pmod 7,\,7\pmod 8,\,8\pmod 9,\,33\pmod {35}\) \\
\bottomrule
\end{longtable}

With these presets, the integer programming filter certifies infeasibility for
the following candidates:
\[
\{
5544,5880,6048,6300,6552,6720,6930,7056,7392,8064,8568,8820,9072,9576
\}.
\]
Thus, after the sum filter and this integer programming filter, the remaining
unresolved candidates below \(10080\) are
\[
    5040,\qquad 7560,\qquad 8400,\qquad 9240.
\]

\section{The Partial Sum Filter}\label{sec:partial-sum-filter}

The third filter is a \emph{partial sum filter}. Like the sum filter, it is based
on reciprocal upper bounds. The difference is that it first computes how much of
\(\Z/L\Z\) can actually be covered by a selected subset of moduli, instead of
assuming that all selected residue classes can be made disjoint. This makes it
strictly more refined than the original sum filter.

Let
\[
    M_L=\{d:d\mid L,\ d\geq 7\}
\]
and split
\[
    M_L=A\sqcup B.
\]
For the subset \(A\), define
\[
    \alpha_A(L)
    :=
    \max_{\{a_m:m\in A\}}
    \frac{1}{L}
    \left|
        \bigcup_{m\in A}
        \{x\in \Z/L\Z:x\equiv a_m\pmod m\}
    \right|.
\]
Thus \(\alpha_A(L)\) is the largest possible proportion of residues modulo
\(L\) that can be covered using only the moduli in \(A\). For the remaining
moduli \(B\), even in the most optimistic case, their contribution is at most
\[
    \sum_{m\in B}\frac1m.
\]
Therefore every full divisor-completed system satisfies
\[
    \frac{1}{L}
    \left|
        \bigcup_{m\in M_L}
        \{x\in \Z/L\Z:x\equiv a_m\pmod m\}
    \right|
    \leq
    \alpha_A(L)+\sum_{m\in B}\frac1m.
\]
Consequently, if
\[
    \alpha_A(L)+\sum_{m\in B}\frac1m<1,
\]
then \(L\) is impossible.

The value \(\alpha_A(L)\) is computed by a smaller integer program. For each
\(m\in A\) and each residue \(a\pmod m\), introduce a binary variable
\(x_{m,a}\), with the same meaning as in Section~\ref{sec:ip-filter}. Impose
\[
    \sum_{a=0}^{m-1}x_{m,a}=1,
    \qquad m\in A.
\]
For each \(b\in\Z/L\Z\), introduce a binary variable \(y_b\), where \(y_b=1\)
means that \(b\) is covered by the selected classes from \(A\). The constraint
\[
    y_b
    \leq
    \sum_{m\in A}x_{m,\,b\bmod m}
\]
forces \(y_b=0\) whenever \(b\) is not covered. Maximizing
\[
    \sum_{b=0}^{L-1}y_b
\]
then gives \(L\alpha_A(L)\). In practice, it is sometimes enough to use the
branch-and-bound upper bound: if
\[
    \frac{\mathrm{BestBound}}{L}+\sum_{m\in B}\frac1m<1,
\]
then infeasibility follows without solving the full covering feasibility problem.

This partial sum filter rules out \(L=8400\) and \(L=9240\). We record the
\(L=9240\) case explicitly. The divisor-completed modulus set is
\[
\begin{aligned}
M_{9240}=\{&
7,8,10,11,12,14,15,20,21,22,24,28,30,33,35,40,42,44,55,56,\\
&60,66,70,77,84,88,105,110,120,132,140,154,165,168,210,220,\\
&231,264,280,308,330,385,420,440,462,616,660,770,840,924,\\
&1155,1320,1540,1848,2310,3080,4620,9240\}.
\end{aligned}
\]
Take the first twenty moduli as
\[
\begin{aligned}
A=\{&
7,8,10,11,12,14,15,20,21,22,24,28,30,33,35,40,42,44,55,56
\},
\end{aligned}
\]
and set \(B=M_{9240}\setminus A\). Then
\[
    \sum_{m\in B}\frac1m
    =
    \frac{877}{4620}
    \approx 0.1898268.
\]
It is therefore enough to prove
\[
    \alpha_A(9240)
    <
    1-\frac{877}{4620}
    =
    \frac{3743}{4620}
    \approx 0.8101732.
\]
The partial-cover optimization certifies this inequality. Hence \(L=9240\)
cannot support a covering system. The same method also rules out \(L=8400\).

After the partial sum filter, only two candidates below \(10080\) remain:
\[
    5040,\qquad 7560.
\]
These two cases are resistant to the partial sum method because they admit very
dense partial covers. The following examples are not covering systems; they only
illustrate why a coarse upper-bound certificate is insufficient.

\begin{example}\label{ex:near-cover-5040}
For \(L=5040\), the following residue classes give a dense partial cover:
\[
\begin{array}{llll}
6\pmod 7, & 7\pmod 8, & 8\pmod 9, & 3\pmod {10},\\
0\pmod {12}, & 8\pmod {14}, & 9\pmod {15}, & 3\pmod {16},\\
2\pmod {18}, & 5\pmod {20}, & 4\pmod {21}, & 6\pmod {24},\\
2\pmod {28}, & 27\pmod {30}, & 0\pmod {35}, & 29\pmod {36},\\
1\pmod {40}, & 10\pmod {42}, & 14\pmod {45}, & 43\pmod {48},\\
26\pmod {56}, & 37\pmod {60}, & 40\pmod {63}, & 56\pmod {70},\\
32\pmod {72}, & 11\pmod {80}, & 16\pmod {84}, & 77\pmod {90},\\
49\pmod {105}, & 10\pmod {112}, & 21\pmod {120}, & 82\pmod {126},\\
112\pmod {140}, & 68\pmod {144}, & 18\pmod {168}, & 109\pmod {180},\\
28\pmod {210}, & 75\pmod {240}, & 124\pmod {252}, & 42\pmod {280}.
\end{array}
\]
It covers \(4857\) residue classes in \(\Z/5040\Z\), leaving \(183\) residue
classes uncovered. Thus it covers about \(96.37\%\) of all residues modulo
\(5040\).
\end{example}

\begin{example}\label{ex:near-cover-7560}
For \(L=7560\), the following residue classes give another dense partial cover:
\[
\begin{array}{llll}
0\pmod {7}, & 0\pmod {8}, & 3\pmod {9}, & 7\pmod {10},\\
2\pmod {12}, & 3\pmod {14}, & 3\pmod {15}, & 4\pmod {18},\\
11\pmod {20}, & 13\pmod {21}, & 20\pmod {24}, & 0\pmod {27},\\
25\pmod {28}, & 29\pmod {30}, & 15\pmod {35}, & 10\pmod {36},\\
1\pmod {40}, & 19\pmod {42}, & 24\pmod {45}, & 36\pmod {54},\\
37\pmod {56}, & 23\pmod {60}, & 16\pmod {63}, & 65\pmod {70},\\
28\pmod {72}, & 67\pmod {84}, & 6\pmod {90}, & 43\pmod {105},\\
18\pmod {108}, & 113\pmod {120}, & 124\pmod {126}, & 60\pmod {135},\\
95\pmod {140}, & 5\pmod {168}, & 42\pmod {180}, & 25\pmod {189},\\
125\pmod {210}, & 180\pmod {216}, & 178\pmod {252}, & 9\pmod {270},\\
141\pmod {280}, & 169\pmod {315}, & 132\pmod {360}, & 88\pmod {378},\\
215\pmod {420}, & 52\pmod {504}, & 99\pmod {540}, & 589\pmod {630},\\
718\pmod {756}, & 293\pmod {840}, & 510\pmod {945}, & 150\pmod {1080},\\
379\pmod {1260}, & 340\pmod {1512}, & 1665\pmod {1890}, & 2301\pmod {2520},\\
2805\pmod {3780}, & 3660\pmod {7560}.
\end{array}
\]
A direct computation modulo \(7560\) shows that these residue classes cover
\(7416\) of the \(7560\) residue classes. Thus the covered proportion is
\[
    \frac{7416}{7560}=\frac{103}{105}\approx 0.98095,
\]
so only \(144\) residue classes remain uncovered.
\end{example}

\section{Final Gurobi Computation}\label{sec:gurobi}

It remains to settle \(L=5040\) and \(L=7560\), and to verify the feasible case
\(L=10080\). We use the full divisor-completed feasibility formulation from
Section~\ref{sec:ip-filter}, with the fixed classes listed there. The
computations were performed with Gurobi~13.0.1 \cite{Gurobi2026} on a multi-core
Linux server.\footnote{The computations were carried out on a Dell PowerEdge
T640 server running Ubuntu 20.04.6 LTS, equipped with two Intel Xeon Gold 6258R
CPUs. At most \(32\) Gurobi threads were used in each run, and the memory usage
remained below approximately \(10\)GB. The code used for these computations is
available at
\url{https://github.com/zhangshiliang502-droid/Distinct-Covering-System-With-m-7}.}

The final complete runs certify infeasibility for \(L=5040\) and \(L=7560\).
For reproducibility, the table below records the full Gurobi outcomes for all
sum-filter candidates, as well as the feasible value \(L=10080\). An infeasible
status means that no choice of residue classes for the corresponding
divisor-completed modulus set gives a complete cover modulo \(L\).

\begin{longtable}{c c c c c}
\toprule
\(L\) & Status & Runtime & Nodes & Simplex iterations \\
\midrule
5040  & infeasible & \(13789.26\)s & \(98459\) & \(530782259\) \\
5544  & infeasible & \(5.34\)s & -- & -- \\
5880  & infeasible & \(7.91\)s & -- & -- \\
6048  & infeasible & \(6.88\)s & -- & -- \\
6300  & infeasible & \(21.05\)s & -- & -- \\
6552  & infeasible & \(6.46\)s & -- & -- \\
6720  & infeasible & \(132.31\)s & -- & -- \\
6930  & infeasible & \(5.13\)s & -- & -- \\
7056  & infeasible & \(5.78\)s & -- & -- \\
7392  & infeasible & \(6.99\)s & -- & -- \\
7560  & infeasible & \(13159.68\)s & \(13533\) & \(186697085\) \\
8064  & infeasible & \(7.24\)s & -- & -- \\
8400  & infeasible & \(922.05\)s & \(43\) & \(1029894\) \\
8568  & infeasible & \(7.08\)s & -- & -- \\
8820  & infeasible & \(17.90\)s & -- & -- \\
9072  & infeasible & \(9.54\)s & -- & -- \\
9240  & infeasible & \(590.53\)s & \(19\) & \(574977\) \\
9576  & infeasible & \(9.03\)s & -- & -- \\
10080 & feasible; independently verified & \(9419.56\)s & \(122934\) & \(212547442\) \\
\bottomrule
\end{longtable}

We now give the explicit covering system at \(L=10080\). One checks modulo
\(10080\) that the following distinct congruence classes cover all integers:
\[
\begin{array}{rrrr}
6 \pmod 7, &
7 \pmod 8, &
8 \pmod 9, &
6 \pmod {10},\\
9 \pmod {12}, &
8 \pmod {14}, &
12 \pmod {15}, &
3 \pmod {16},\\
14 \pmod {18}, &
0 \pmod {20}, &
4 \pmod {21}, &
13 \pmod {24},\\
26 \pmod {28}, &
24 \pmod {30}, &
27 \pmod {32}, &
33 \pmod {35},\\
5 \pmod {36}, &
11 \pmod {40}, &
16 \pmod {42}, &
2 \pmod {45},\\
1 \pmod {48}, &
52 \pmod {56}, &
30 \pmod {60}, &
7 \pmod {63},\\
28 \pmod {70}, &
29 \pmod {72}, &
59 \pmod {80}, &
10 \pmod {84},\\
74 \pmod {90}, &
43 \pmod {96}, &
93 \pmod {105}, &
73 \pmod {112},\\
64 \pmod {120}, &
112 \pmod {126}, &
38 \pmod {140}, &
65 \pmod {144},\\
43 \pmod {160}, &
121 \pmod {168}, &
110 \pmod {180}, &
18 \pmod {210},\\
169 \pmod {224}, &
185 \pmod {240}, &
154 \pmod {252}, &
248 \pmod {280},\\
203 \pmod {288}, &
128 \pmod {315}, &
313 \pmod {336}, &
209 \pmod {360},\\
292 \pmod {420}, &
75 \pmod {480}, &
217 \pmod {504}, &
404 \pmod {560},\\
578 \pmod {630}, &
505 \pmod {672}, &
281 \pmod {720}, &
472 \pmod {840},\\
553 \pmod {1008}, &
233 \pmod {1120}, &
532 \pmod {1260}, &
875 \pmod {1440},\\
124 \pmod {1680}, &
281 \pmod {2016}, &
2044 \pmod {2520}, &
2153 \pmod {3360},\\
5033 \pmod {5040}, &
7193 \pmod {10080}.
\end{array}
\]
The moduli are distinct, the least modulus is \(7\), and the least common
multiple is \(10080\). Hence \(L_{\min}(7)\leq 10080\).

On the other hand, Sections~\ref{sec:sum-filter}--\ref{sec:partial-sum-filter}
and the complete computations above rule out every possible \(L<10080\).
Therefore \(L_{\min}(7)\geq 10080\). Combining the two inequalities proves
Theorem~\ref{thm:main}.

\section{Conclusion}\label{sec:conclusion}

We have constructed a distinct covering system with minimum modulus \(7\) and
least common multiple \(10080\), improving the previously conjectured value
\(15120\). The minimality proof is obtained by a successive filtering strategy:
the reciprocal-sum filter reduces the raw multiples of \(7\) to eighteen
candidates; the integer-programming filter eliminates most of them; the
partial-sum filter removes the difficult cases \(8400\) and \(9240\); and final
Gurobi computations rule out \(5040\) and \(7560\). Thus the optimal value is
\[
    L_{\min}(7)=10080.
\]

The partial-sum filter also suggests a useful direction for future work. The
near-cover examples for \(5040\) and \(7560\) indicate that some candidates are
difficult because they can cover almost all residue classes without forming a
complete covering system. Sharper hybrid certificates for separating such
near-covers from genuine coverings may reduce the need for full branch-and-bound
verification in later cases.

\newpage

\section*{Acknowledgements}

This work was supported by ShanghaiTech AI Initiative
(Grant No. AI2026B16). We thank Hao Chen, Lian Duan, Boqing Xue, and Jiayu
Zhai for helpful support, discussions, and encouragement.

\ContactInfo

\begin{thebibliography}{99}
\bibitem{Erdos1950}
P. Erd\H{o}s,
\emph{On integers of the form $2^k+p$ and some related problems},
\emph{Summa Brasil. Math.} \textbf{2} (1950), 113--123.

\bibitem{Swift1954}
J. D. Swift,
\emph{Sets of covering congruences},
\emph{Bull. Amer. Math. Soc.} \textbf{60} (1954), 390.

\bibitem{Churchhouse1968}
R. F. Churchhouse,
\emph{Covering sets and systems of congruences},
in \emph{Computers in Mathematical Research},
North-Holland, Amsterdam, 1968, pp.~20--36.

\bibitem{Choi1971}
S. L. G. Choi,
\emph{Covering the set of integers by congruence classes of distinct moduli},
\emph{Math. Comp.} \textbf{25} (1971), 885--895.
\href{https://doi.org/10.1090/S0025-5718-1971-0297692-7}
{doi:10.1090/S0025-5718-1971-0297692-7}.

\bibitem{Krukenberg1971}
C. E. Krukenberg,
\emph{Covering Sets of the Integers},
Ph.D. thesis, University of Illinois at Urbana--Champaign, 1971.

\bibitem{Morikawa1981}
R. Morikawa,
\emph{On a method to construct covering sets},
\emph{Bull. Fac. Liberal Arts Nagasaki Univ.} \textbf{22} (1981), no.~1, 1--11.

\bibitem{PorubskySchonheim2002}
\v{S}. Porubsk\'y and J. Sch\"onheim,
\emph{Covering systems of Paul Erd\H{o}s: past, present and future},
in \emph{Paul Erd\H{o}s and His Mathematics, I},
J\'anos Bolyai Math. Soc., Budapest, 2002, pp.~581--627.

\bibitem{Gibson2009}
D. J. Gibson,
\emph{A covering system with least modulus $25$},
\emph{Math. Comp.} \textbf{78} (2009), no.~266, 1127--1146.

\bibitem{Nielsen2009}
P. P. Nielsen,
\emph{A covering system whose smallest modulus is $40$},
\emph{J. Number Theory} \textbf{129} (2009), no.~3, 640--666.

\bibitem{Owens2014}
T. Owens,
\emph{A covering system with minimum modulus $42$},
Master's thesis, Wake Forest University, 2014.

\bibitem{Hough2015}
R. Hough,
\emph{Solution of the minimum modulus problem for covering systems},
\emph{Ann. of Math.} \textbf{181} (2015), no.~1, 361--382.

\bibitem{DaltonTrifonov2022}
J. Dalton and O. Trifonov,
\emph{Extreme covering systems},
\emph{J. Integer Seq.} \textbf{25} (2022), no.~9,
Article~22.9.1, 28 pp.

\bibitem{Balister2024}
P. Balister,
\emph{Erd\H{o}s covering systems},
in \emph{Surveys in Combinatorics 2024},
Cambridge University Press, Cambridge, 2024, pp.~31--54.

\bibitem{Klein2026}
J. Klein,
\emph{On a conjecture of Krukenberg and a problem of Dalton and Trifonov},
\emph{Integers} \textbf{26} (2026), \#A38.
\href{https://doi.org/10.5281/zenodo.19402698}
{doi:10.5281/zenodo.19402698}.

\bibitem{Gurobi2026}
Gurobi Optimization, LLC,
\emph{Gurobi Optimizer Reference Manual},
version 13.0.1, 2026.
\url{https://www.gurobi.com}.

\end{thebibliography}
\end{document}